\documentclass{amsart}

\usepackage{graphicx,mathrsfs,amssymb,amsfonts}

\usepackage[pdfdisplaydoctitle=true,colorlinks=true,urlcolor=blue,citecolor=blue,%
linkcolor=blue,pdfstartview=FitH,pdfpagemode=None,bookmarksnumbered=true]{hyperref}
\usepackage{cmap} % make the PDF file "searchable and copyable"
\usepackage[usenames]{color}
\usepackage{leftidx}
\usepackage[normalem]{ulem}

\numberwithin{equation}{section}

\newtheorem{thm}{Theorem}[section]
\newtheorem{cor}[thm]{Corollary}
\newtheorem{lem}[thm]{Lemma}
\newtheorem{prop}[thm]{Proposition}
\theoremstyle{definition}
\newtheorem{defn}[thm]{Definition}
\theoremstyle{remark}
\newtheorem{rem}[thm]{Remark}

\numberwithin{equation}{section}

\newcommand{\cz}{{\mathbb C}}
\newcommand{\nz}{{\mathbb N}}
\newcommand{\rz}{{\mathbb R}}

\newcommand{\calA}{\mathcal{A}}
\newcommand{\calE}{\mathcal{E}}
\newcommand{\calP}{\mathcal{P}}
\newcommand{\calQ}{\mathcal{Q}}

\newcommand{\scrC}{\mathscr{C}}
\newcommand{\scrF}{\mathscr{F}}
\newcommand{\scrL}{\mathscr{L}}
\newcommand{\scrS}{\mathscr{S}}

\newcommand{\cl}{\mathrm{cl}}
\newcommand{\dbar}{d\hspace*{-0.08em}\bar{}\hspace*{0.1em}}

\newcommand{\lra}{\longrightarrow}
\newcommand{\op}{\mathrm{op}}
\newcommand{\spk}[1]{\langle#1\rangle}
\newcommand{\wh}{\widehat}
\newcommand{\wt}{\widetilde}

\begin{document}

\title[Fredholm property of bisingular operators]{On the Fredholm property of \\ 
bisingular pseudodifferential operators}

%----------Author 1
\author{Massimo Borsero}
\address{Universit\`a degli Studi di Torino\\ Dipartimento di Matematica ``Giuseppe Peano''\\ 
Via Carlo Alberto 10 \\ 10123, Torino (Italy)}
\email{massimo.borsero@unito.it}

%----------Author 2
\author{J\"org Seiler}
\address{Universit\`a degli Studi di Torino\\ Dipartimento di Matematica ``Giuseppe Peano''\\ 
Via Carlo Alberto 10 \\ 10123, Torino (Italy)}
\email{joerg.seiler@unito.it}



%\date{January 1, 2004}

\begin{abstract}
For operators belonging either to a class of global bisingular pseudodifferential operators on $\rz^m\times\rz^n$ 
or to a class of bisingular pseudodifferential operators on a product $M\times N$ of two closed smooth manifolds, 
we show the equivalence of their ellipticity (defined by the invertibility of certain operator-valued,  
homogeneous principal symbols) and their Fredholm mapping property in associated scales of Sobolev spaces. 
We also prove the spectral invariance of these operator classes and then extend these results to larger  
classes of Toeplitz type operators. 
\end{abstract}

\maketitle

%%%%%%%%%%%%%%%%%%%%%%%%%%%%%%%%%%%%%%%%%%%%%%%%%%
%%%%%%%%%%%%%%%%%%%%%%%%%%%%%%%%%%%%%%%%%%%%%%%%%%
\section{Introduction}\label{sec:intro}

Calculi of bisingular pseudodifferential operators can be seen as a systematic approach for studying tensor products 
of pseudodifferential operators. Focusing on elliptic theory, a typical question would be the following: Given 
classical $($or poly-homogeneous$)$ pseudodifferential operators $A_j\in L^\mu_\cl(M)$ and $B_j\in L^\nu_\cl(N)$ 
for $j=1,\ldots,k$, on smooth manifolds $M$ and $N$, how can we characterize the existence of a parametrix, 
the Fredholm property or the invertibilty of the operator $A_1\otimes B_1+\ldots+A_k\otimes B_k$? Here, the   
tensor product $A\otimes B$ denotes an operator acting on functions defined on $M\times N$ with the property 
that 
 $$ A\otimes B(u\otimes v)=Au\otimes Bv,\qquad u\in\scrC^\infty(M),\;v\in\scrC^\infty(N),$$
where $(f\otimes g)(x,y)=f(x)g(y)$ for any two functions $f$ and $g$ on $M$ and $N$, respectively.  
Such tensor products, in general, do not define a classical pseudodifferential operator on $M\times N$, hence the 
question cannot be answered using only the standard calculus. 

Questions of this kind are not only of academic interest but arose, in particular, naturally in the framework of the 
famous Atiyah-Singer index theorem. In fact, Atiyah and Singer in \cite{AtSi} were led to study systems of the form 
 $$A\boxtimes B=\begin{pmatrix}A\otimes 1 & -1\otimes B^*\\1\otimes B & A^*\otimes 1\end{pmatrix},$$
where both $A$ and $B$ are zero-order classical 
pseudodifferential operators on $M$ and $N$, respectively. Again, $A\boxtimes B$ is not a classical pseudodifferential 
operator on $M\times N$. However, if both $A$ and $B$ are elliptic, then $A\boxtimes B$ is a Fredholm operator in 
$L^2(M\times N,\cz^2)$ with index $\mathrm{ind}\,A\boxtimes B=\mathrm{ind}\,A\cdot\mathrm{ind}\,B$.  

Motivated by these phenomena, Rodino in \cite{Rodi} introduced a pseudodifferential calculus of operators 
acting on sections of vector bundles over a product of smooth, closed $($i.e., compact and without boundary$)$ 
manifolds $M\times N$, containing such kinds of tensor product type operators. 
We recall the main features and ideas in Section \ref{sec:03}. In this calculus, operators 
can be composed and parametrices to elliptic elements can be constructed. Ellipticity in this context refers to the 
invertibility of two \emph{operator-valued} principal symbols associated with each operator $($roughly speaking, 
each such principal symbol is defined on the co-tangent bundle of one of the two manifolds and takes values in 
the space of classical pseudodifferential operators of the other manifold$)$. 
In Section \ref{sec:03.1.2} we carefully discuss these principal symbols, 
developing a formalism necessary for our application to so-called Toeplitz type operators presented in 
Section \ref{sec:04}. 

As a consequence of the existence of parametrices to elliptic operators, as shown in \cite{Rodi}, 
elliptic operators act as 
Fredholm operators in a certain scale of naturally associated $L^2$-Sobolev spaces. 
The main result in the present 
paper is the proof of the reverse statement: 
If a bisingular pseudodifferential operator in the calculus of \cite{Rodi} is Fredholm
it necessarily must be elliptic. In other words, the ellipticity condition used in the calculus is ``optimal". 
The method of our proof is based on techniques introduced in Gohberg \cite{Gohb} and H\"ormander \cite{Horm}. 
Also, as a consequence, we obtain that the calculus of Rodino is spectrally invariant.
{Both equivalence of Fredholm property and ellipticity as well as the spectral invariance have been employed in 
the very recent work Bohlen \cite{Bohl}, where the meromorphic structure of the $\eta$-function for 
$($scaler-valued$)$ bisingular pseudodifferential operators is investigated}.  

Of course one can pose analogous questions also in case where $M$ and $N$ are not compact. It then depends very 
much on the sort of non-compactness which kind of operators one would consider. In the present paper, we 
investigate the case $M=\rz^m$ and $N=\rz^n$ and work with bisingular operators based on pseudodifferential 
operators of Shubin type, cf.\ \cite{Shub}. Such a calculus was recently considered in 
Battisti, Gramchev, Rodino and Pilipovi\'{c} \cite{BGRP}, where a Weyl law for the spectral counting
function of global bisingular operators has been obtained, and also in Nicola and Rodino \cite{NR06}, where the 
noncommutative residue is studied. Again we show, in Section \ref{sec:02}, 
equivalence of ellipticity and Fredholm property as well as spectral invariance. 

As a matter of fact, our results allow us to treat even more general kinds of bisingular operators, of so-called 
\emph{Toeplitz type}, both in the context of bisingular operators on $M\times N$ and $\rz^m\times\rz^n$, 
repectively. To this end we show in Section \ref{sec:04} that general results of Seiler \cite{Seil} on abstract 
pseudodifferential operators of Toeplitz type apply in the present two settings of bisingular operator classes.
{As an application, we prove the existence of bisingular order-reductions.}

{
The addressed question of \emph{characterizing} the Fredholm property of pseudodifferential operators in terms of 
the invertibility of associated principal symbols is a fundamental problem whenever working with algebras/calculi of 
pseudodifferential operators. In many concrete cases such results are valid; we just mention the calculi of Schulze 
\cite{Schu} for manifolds $($with and without boundary$)$ with conical singularities, edges, and higher singularities, 
and the calculi of Melrose \cite{Melr} for corner manifolds. A general approach to this question, which contains  
many of these calculi as specific examples,  has been developed by Nistor and co-authors in the framework of 
pseudodifferential operators on groupoids, see \cite{LMN} and references therein. In \cite{Mant}, Mantoiu uses 
$C^*$-algebra techniques to investigate the essential spectrum $($Fredholm spectrum$)$ of Schr\"odinger operators 
on locally compact Lie groups, including bisingular Schr\"odinger operators as particular examples. 

Given a specific pseudodifferential calculus, one may be interested in a corresponding calculus 
of bisingular operators and study the relation between ellipticity and Fredholm property. In this perspective, our 
paper only concerns a relatively simple situation; more complicated settings might be subject to future research.     
} 

%\textbf{Acknowledgements.} 
%We thank the anonymous referee for valuable remarks and suggestions that led to various improvements of the paper.  

%%%%%%%%%%%%%%%%%%%%%%%%%%%%%%%%%%%%%%%%%%%%%%%%%%
%%%%%%%%%%%%%%%%%%%%%%%%%%%%%%%%%%%%%%%%%%%%%%%%%%
\section{Bisingular operators of Shubin type}\label{sec:02}

In the present section we show the equivalence of ellipticity and Fredholm property for a certain class of 
global bisingular operators on $\rz^{m}\times\rz^{n}$, a bisingular version of operators of Shubin type \cite{Shub}. 
For the more technical details of this calculus we refer the reader to the recent paper 
\cite{BGRP}.\footnote{Actually, in \cite{BGRP} the authors work with a class of symbols slightly larger than 
the one employed here. They only require the existence of the homogeneous principal symbols while we ask the 
existence of  complete asymptotic expansions in homogeneous components. However, our approach carries over 
without modification to this larger calculus  and our results, i.e., Theorems \ref{thm:shubin}, \ref{thm:spectral} and 
Corollary \ref{cor:extension}, remain valid. In fact, our calculus coincides with the one of \cite{NR06}, where it is 
presented with a slightly different formalism. \label{fn:01}} 

{
Let us introduce here two notations which we will use throughout the whole paper. We write 
$\spk{y}=(1+|y|^2)^{1/2}$ for vectors $y\in\rz^k$. 
In case $y=(y_1,y_2)$ we shall also write $\spk{y_1,y_2}:=\spk{(y_1,y_2)}$. 

Moreover, the unit-sphere in $\rz^k$ we shall denote by $\mathbb{S}^{k-1}$. 
}

%%%%%%%%%%%%%%%%%%%%%%%%%%%%%%%%%%%%%%%%%%%%%%%%%%
\subsection{Shubin type symbols with values in a Fr\'{e}chet space}\label{sec:02.1}

Let $F$ be a Fr\'{e}chet space with topology given by the system of semi-norms $p_0,p_1,p_2,\ldots$. 

For $\nu\in\rz$ we let $\Gamma^\nu(\rz^n;F)$ denote the space of all smooth functions $a:\rz^n\times\rz^n\to F$ 
satisfying, for any $k\in\nz$,   
\begin{align}\label{eq:symbol1} 
 q_k(a):=\sup_{\substack{x,\xi\in\rz^n\\ j+|\alpha|+\beta|\le k}}
 p_j\big(D^\alpha_\xi D^\beta_x a(x,\xi)\big)\spk{x,\xi}^{|\alpha|+|\beta|-\nu}<+\infty.
\end{align}
These semi-norms turn $\Gamma^\nu(\rz^n;F)$ into a Fr\'{e}chet space. 

The subspace $\Gamma^\nu_\cl(\rz^n;F)$ of classical $($or poly-homogeneous$)$ symbols consists of those 
elements of $\Gamma^\nu(\rz^n;F)$ for which there exist smooth functions 
\begin{equation}\label{eq:homogeneous}
 a^{(\nu-j)}:(\rz^n\times\rz^n)\setminus \{0\}\to F, \qquad j=0,1,2,\ldots,
\end{equation}
that are positively homogeneous of degree $\nu-j$ in $(x,\xi)$, i.e., 
 $$a^{(\nu-j)}(tx,t\xi)=t^{\nu-j}\,a^{(\nu-j)}(x,\xi)\qquad \forall\;t>0\quad\forall\;(x,\xi)\not=0,$$
such that  
 $$r_N(a):=a-\sum_{j=0}^{N-1}\chi a^{(\nu-j)}\;\in\;\Gamma^{\nu-N}(\rz^n;F)
     \qquad \forall\; N=0,1,2,\ldots,$$
where $\chi(x,\xi)$ is a smooth zero-excision function, i.e., $\chi\equiv0$ near the origin and $1-\chi$ has compact 
support. Note that the \emph{homogeneous components} $a^{(\nu-j)}$ are uniquely determined by $a$; 
the component 
$a^{(\nu)}$ is called the \emph{homogeneous principal symbol} of $a$. By homogeneity, we may identify every  
component with a smooth, $F$-valued function defined on the unit-sphere $\mathbb{S}^{2n-1}$ in 
$\rz^n\times\rz^n$. Then the maps 
\begin{align*}
 a\mapsto r_N(a)&:\Gamma^{\nu}_\cl(\rz^n;F)\lra \Gamma^{\nu-N}(\rz^n;F), \\
 a\mapsto a^{(\nu-j)}&:\Gamma^{\nu}_\cl(\rz^n;F)\lra \scrC^\infty(\mathbb{S}^{2n-1};F) 
\end{align*}
with $j,N=0,1,2,\ldots$, induce a Fr\'echet topology on $\Gamma^{\nu}_\cl(\rz^n;F)$. 

Finally, note that 
 $$\Gamma^{-\infty}(\rz^n;F):=\mathop{\mbox{\Large$\cap$}}_{\nu\in\rz}\Gamma^{\nu}(\rz^n;F)
     =\mathop{\mbox{\Large$\cap$}}_{\nu\in\rz}\Gamma^{\nu}_\cl(\rz^n;F)$$
coincides with the Schwartz space $\scrS(\rz^n,F)$ of rapidly decreasing, $F$-valued functions. 

%%%%%%%%%%%%%%%%%%%%%%%%%%%%%%%%%%%%%%%%%%%%%%%%%%
\subsubsection{Operator-valued symbols}\label{sec:02.1.1}

Of particular importance is the case $F=\scrL(E_1,E_2)$, the Banach space of all bounded, linear operators 
$E_1\to E_2$ between two Hilbert spaces. In this case we associate with $a\in \Gamma^\nu(\rz^n,\scrL(E_1,E_2))$ 
the pseudodifferential operator $A=\op(a):\scrS(\rz^n,E_1)\to \scrS(\rz^n,E_2)$ defined by 
 $$(A u)(x)=\int e^{ix\xi}a(x,\xi)\wh{u}(\xi)\,\dbar\xi,\qquad \scrS(\rz^n,E_1).$$ 
For $E_1=E_2=\cz$ these are the standard pseudodifferential symbols (respectively operators) 
from the Shubin class as introduced in \cite{Shub}. Note that operators associated with symbols of order $-\infty$ 
are integral operators with integral kernels that are Schwartz functions in both variables.  

%%%%%%%%%%%%%%%%%%%%%%%%%%%%%%%%%%%%%%%%%%%%%%%%%%
{
\subsubsection{Ellipticity}\label{sec:02.1.2}

$a\in \Gamma^\nu_\cl(\rz^n,\scrL(E_1,E_2))$ is called \textit{elliptic}, if its homogeneous principal symbol 
$a^{(\nu)}$ from 
\eqref{eq:homogeneous} is invertible for every $(x,\xi)\not=0$. In this case $a$ admits a so-called parametrix, i.e.,  
a symbol $b\in \Gamma^{-\nu}_\cl(\rz^n,\scrL(E_2,E_1))$ such that $\op(a)\op(b)=1-\op(r_1)$ and 
$\op(b)\op(a)=1-\op(r_2)$ with symbols $r_1$ and $r_2$ or order $-\infty$. 
}

%%%%%%%%%%%%%%%%%%%%%%%%%%%%%%%%%%%%%%%%%%%%%%%%%%
{
\subsubsection{Parameter-dependent operators and order-reductions}\label{sec:02.1.3}

In the definition of the symbol classes from the beginning of Section \ref{sec:02.1} one may replace the covariable 
$\xi$ with $\eta:=(\xi,\sigma)$, where $\sigma$ is a real parameter. This then leads to symbol classes 
denoted by $\Gamma^\nu_{(\cl)}(\rz^n,\rz_\sigma;F)$ and to corresponding operator-families $A(\sigma)$ 
in case $F=\scrL(E_1,E_2)$. Ellipticity asks the invertibility of the homogeneous principal symbol for all 
$(x,\eta)\not=0$ and implies the existence of a parameter-dependent parametrix, i.e.,  
$\op(a)(\sigma)\op(b)(\sigma)=1-\op(r_1)(\sigma)$ and 
$\op(b)(\sigma)\op(a)(\sigma)=1-\op(r_2)(\sigma)$ with parameter-dependent $r_1$ and $r_2$ of order $-\infty$. 
Employing that the parameter in $r_1$ and $r_2$ is rapidly decreasing as it tends to $\pm\infty$, one can modify $b$ in such a 
way, that $\op(a)(\sigma)\op(b)(\sigma)-1$ and $\op(b)(\sigma)\op(a)(\sigma)-1$ are compactly supported in $\sigma$. 
In other words, if $a(\sigma)\in \Gamma^\nu_{(\cl)}(\rz^n,\rz_\sigma;\scrL(E_1,E_2))$ is parameter-elliptic and 
$\sigma_0$ is sufficiently large, then 
 $$\lambda^\nu(x,\xi):=a(x,\xi,\sigma_0)\in \Gamma^\nu_{\cl}(\rz^n;\scrL(E_1,E_2))$$
and 
 $$\lambda^{-\nu}(x,\xi):=b(x,\xi,\sigma_0)\in \Gamma^{-\nu}_{\cl}(\rz^n;\scrL(E_2,E_1))$$
satisfy $\op(\lambda^\nu)\op(\lambda^{-\nu})=\mathrm{id}_{E_2}$ and 
$\op(\lambda^{-\nu})\op(\lambda^{\nu})=\mathrm{id}_{E_1}$. Any such $\lambda^\nu$ is called an order-reduction 
of order $\nu$. For example, in case $E=E_1=E_2$ one can take
\begin{equation}\label{eq:redord} 
a(x,\xi,\sigma)=[x,\xi,\sigma]^\nu\,\mathrm{id}_E,
\end{equation} 
where $[\cdot]:\rz^{2n+1}_{x,\xi,\sigma}\to\rz$ denotes a positive smooth function that coincides with the usual 
modulus outside some neighborhood of the origin. 
}

%%%%%%%%%%%%%%%%%%%%%%%%%%%%%%%%%%%%%%%%%%%%%%%%%%
{
\subsubsection{Sobolev spaces}\label{sec:02.1.4}

Let $E$ be a Hilbert space and $\Lambda^{s}=\op(\lambda^s)$ be an order-reduction of order $s$ as described in 
the previous subsection $($with $E=E_0=E_1)$. The Sobolev space $Q^s(\rz^n,E)$ of order $s$ is defined as the 
closure of $\scrS(\rz^n,E)$ with respect to the norm $\|u\|_s=\|\Lambda^su\|_{L^2(\rz^n,E)}$. 

For a symbol $a\in \Gamma^\nu(\rz^n,\scrL(E_1,E_2))$, the associated operator $A=\op(a)$ extends by continuity  
to $A:Q^s(\rz^n,E_1)\to Q^{s-\nu}(\rz^n,E_2)$ for every $s\in\rz$. 
}

%%%%%%%%%%%%%%%%%%%%%%%%%%%%%%%%%%%%%%%%%%%%%%%%%%
\subsection{Bisingular symbols and their calculus}\label{sec:02.2}

Let us denote by 
 $$\Gamma^{\mu,\nu}(\rz^m\times\rz^n;\cz^k,\cz^\ell),\qquad \mu,\nu\in\rz\cup\{-\infty\},\quad k,l\in\nz,$$ 
the space of all smooth functions 
$a:\rz^m\times\rz^m\times\rz^n\times \rz^n\to\cz^{\ell\times k}$ $($taking values in the complex 
$\ell\times k$-matrices, identified with $\scrL(\cz^k,\cz^\ell)$ by using the standard basis of $\cz^k$ and 
$\cz^\ell$, respectively$)$ such that 
\begin{align*}
 (x,\xi)\mapsto a_1(x,\xi):=  \Big( (y,\eta)\mapsto a(x,\xi,y,\eta) \Big)
\end{align*}
defines a Fr\'{e}chet space valued symbol  
\begin{align}\label{eq:symba1}
 a_1\in \Gamma^\mu\big(\rz^m;\Gamma^\nu(\rz^n;\cz^{\ell\times k})\big). 
\end{align}
In this case, 
\begin{align*}
 (y,\eta)\mapsto a_2(y,\eta):=  \Big( (x,\xi)\mapsto a(x,\xi,y,\eta) \Big)  
\end{align*}
defines a symbol
\begin{align}\label{eq:symba2}
  a_2\in \Gamma^\nu\big(\rz^n;\Gamma^\mu(\rz^m;\cz^{\ell\times k})\big). 
\end{align}

\begin{rem}
A function $a$ belongs to $\Gamma^{\mu,\nu}(\rz^m\times\rz^n;\cz^k,\cz^\ell)$ if, and only if, it satisfies 
the uniform estimates 
 $$\|D^\alpha_\xi D^\beta_x D^\gamma_\eta D^\delta_y a(x,\xi,y,\eta)\|_{\cz^{\ell\times k}}\le 
    C_{\alpha\beta}
     \spk{x,\xi}^{\mu-|\alpha|-|\beta|}\spk{y,\eta}^{\nu-|\gamma|-|\delta|}$$
for every order of derivatives. 
\end{rem} 

The spaces of \emph{classical} symbols $\Gamma^{\mu,\nu}_\cl(\rz^m\times\rz^n;\cz^k,\cz^\ell)$ are defined as 
above, replacing $\Gamma^\mu$ and $\Gamma^\nu$ by $\Gamma^\mu_\cl$ and $\Gamma^\nu_\cl$, respectively. 

%%%%%%%%%%%%%%%%%%%%%%%%%%%%%%%%%%%%%%%%%%%%%%%%%%
\subsubsection{Operators and Sobolev spaces}\label{sec:02.2.1}

With $a\in \Gamma^{\mu,\nu}(\rz^m\times\rz^n;\cz^k,\cz^\ell)$ we associate, as usual, 
its pseudodifferential operator 
\begin{align}\label{eq:op_schwarz}
 A=\op(a):\scrS(\rz^m\times\rz^n,\cz^k)\lra \scrS(\rz^m\times\rz^n,\cz^\ell).
\end{align}
The map $a\mapsto\op(a)$ establishes a bijection between the respective spaces of symbols and operators. 
Therefore we shall not introduce a new notation for the spaces of operators, but simply write 
$A\in \Gamma^{\mu,\nu}(\rz^m\times\rz^n;\cz^k,\cz^\ell)$. 
Operators of order $(-\infty,-\infty)$ we shall refer to as \emph{regularizing} or \emph{smoothing} operators. 

\begin{rem}
With $A=\op(a)\in\Gamma^\nu(\rz^n)$ and $B=\op(b)\in\Gamma^\mu(\rz^m)$, let 
$a\otimes b\in\Gamma^{\mu,\nu}(\rz^m\times\rz^n)$ be defined by $a\otimes b(x,\xi,y,\eta)=a(x,\xi)b(y,\eta)$. 
The associated operator shall be denoted by $A\otimes B=\op(a\otimes b)$. 
If $u(x,y)=v(x)w(y)$ with rapidly decreasing functions $v$ and $w$, then 
 $$[(A\otimes B) u](x,y)=(Av)(x)(Bw)(y).$$
Such tensor-products, respectively finite linear combinations, are the simplest examples of bisingular operators. 
Using the nuclearity of $\Gamma^\nu_\cl(\rz^n)$ indeed it can be shown that 
\begin{align}\label{eq:tensor}
 \Gamma^{\mu,\nu}_\cl(\rz^m\times\rz^n)=\Gamma^\mu_\cl(\rz^m)\,\wh{\otimes}_\pi \Gamma^\nu_\cl(\rz^n),
\end{align}
where $E\,\wh{\otimes}_\pi F$ denotes the completed, projective tensor-product of two Fr\'{e}chet spaces 
$E$ and $F$, cf.\ \cite{Tre67}.  Note that an equality as in \eqref{eq:tensor} does not hold for the spaces of 
non-classical symbols. 
\end{rem} 

The operator from \eqref{eq:op_schwarz} extends continuously to 
\begin{align}\label{eq:op_sobolev}
 A:Q^{s,t}(\rz^m\times\rz^n,\cz^k)\lra Q^{s-\mu,t-\nu}(\rz^m\times\rz^n,\cz^\ell),\qquad s,t\in\rz,
\end{align}
where $Q^{s,t}(\rz^m\times\rz^n,\cz^j)$ is the $j$-fold sum of $Q^{s,t}(\rz^m\times\rz^n)$, the latter being 
the closure of $\scrS(\rz^m\times\rz^n)$ with respect to the norm 
$u\mapsto\|\Lambda^{s,t}u\|_{L^2(\rz^m\times\rz^n)}$, where  
{
$\Lambda^{s,t}=\Lambda^s_m\otimes\Lambda^t_n$  with order-reductions $\Lambda^s_m$ and $\Lambda^t_n$ 
of order $s$ and $t$ on $\rz^m$ and $\rz^n$, respectively, as described in Section \ref{sec:02.1.3}.
}

Bisingular symbols behave well under composition and taking the formal adjoint, in the sense that$:$
\begin{itemize}
 \item[(1)] Composition of operators, $(A_2,A_1)\mapsto A_2A_1$, induces maps 
   \begin{align*}
    \Gamma^{\mu_2,\nu_2}&(\rz^m\times\rz^n;\cz^j,\cz^\ell)\times 
    \Gamma^{\mu_1,\nu_1}(\rz^m\times\rz^n;\cz^k,\cz^j)\\
    &\lra \Gamma^{\mu_1+\mu_2,\nu_1+\nu_2}(\rz^m\times\rz^n;\cz^k,\cz^\ell).
   \end{align*}
 \item[(2)] Taking the formal $L^2$-adjoint, $A\mapsto A^*$, induces maps 
   \begin{align*}
    \Gamma^{\mu,\nu}&(\rz^m\times\rz^n;\cz^k,\cz^\ell)\lra 
    \Gamma^{\mu,\nu}(\rz^m\times\rz^n;\cz^\ell,\cz^k). 
   \end{align*}
\end{itemize}
The analogous statements are true for classical symbols. 

%%%%%%%%%%%%%%%%%%%%%%%%%%%%%%%%%%%%%%%%%%%%%%%%%%
\subsubsection{Classical symbols and ellipticity}\label{sec:02.2.2}

With a classical operator $A=\op(a)$ belonging to $\Gamma^{\mu,\nu}_\cl(\rz^m\times\rz^n;\cz^k,\cz^\ell)$ 
we associate two principal symbols 
\begin{align*}
  \sigma_1^\mu(A)&=a^{(\mu)}_1\in 
    \scrC^\infty\big(\mathbb{S}^{2m-1},\Gamma^\nu_\cl(\rz^n;\cz^{\ell\times k})\big), \\ 
  \sigma_2^\nu(A)&=a^{(\nu)}_2\in 
    \scrC^\infty\big(\mathbb{S}^{2n-1},\Gamma^\mu_\cl(\rz^m;\cz^{\ell\times k})\big),  
\end{align*}
the homogeneous principal symbol of $a_1$ and $a_2$ as defined in \eqref{eq:symba1} and \eqref{eq:symba2}, 
respectively, restricted to the corresponding unit-sphere. 
Note that 
\begin{align*}%\label{eq:qs}
    \sigma_1^\mu(A)\in 
    \scrC^\infty\big(\mathbb{S}^{2m-1},\scrL(Q^s(\rz^m,\cz^k),Q^{s-\mu}(\rz^m,\cz^\ell))\big),\qquad s\in\rz,
\end{align*}
and similarly for $\sigma_2^\nu(A)$. For compostion and adjoints of operators we have, 
using notation from $(1)$ and $(2)$ above, 
 $$\sigma_1^{\mu_1+\mu_2}(A_2A_1)=\sigma_1^{\mu_2}(A_2)\sigma_1^{\mu_1}(A_1),
     \qquad\sigma_1^{\mu}(A^*)=\sigma_1^{\mu}(A)^*,$$
where the $*$ on the right-hand side is the  formal $L^2$-adjoint 
$\Gamma^\nu(\rz^n;\cz^k,\cz^\ell)\to \Gamma^\nu(\rz^n;\cz^\ell,\cz^k)$. 
Analogous equations hold for the other principal symbol $\sigma_2$. 

\begin{defn}\label{def:ell}
$A\in \Gamma^{\mu,\nu}_\cl(\rz^m\times\rz^n;\cz^k,\cz^k)$ is called elliptic if both $\sigma_1^\mu(A)$ and 
$\sigma_2^\nu(A)$ take values in the invertible operators. 
\end{defn}

In the previous definition, invertibility of $\sigma_1^\mu(A)(x,\xi)$ refers either to invertibility in 
$\scrL(Q^s(\rz^m,\cz^k),Q^{s-\mu}(\rz^m,\cz^k))$ 
for some $s\in\rz$ or to invertibilty in $\Gamma^\nu_\cl(\rz^n;\cz^{k\times k})$, 
{i.e., having an inverse belonging to $\Gamma^{-\nu}_\cl(\rz^n;\cz^{k\times k})$}. 
Due to the spectral invariance of the standard Shubin class $($which is a particular 
case of the spectral invariance of bisingular operators that we shall prove in this paper$)$ both 
possibilities are equivalent. 

The following theorem is one of the main results for elliptic operators: 

\begin{thm}\label{thm:parametrix}
An operator $A\in \Gamma^{\mu,\nu}_\cl(\rz^m\times\rz^n;\cz^k,\cz^k)$ is elliptic if, and only if, there exists 
an operator $B\in \Gamma^{-\mu,-\nu}_\cl(\rz^m\times\rz^n;\cz^k,\cz^k)$ such that 
 $$1-AB, 1-BA \;\in\; \Gamma^{-\infty,-\infty}(\rz^m\times\rz^n;\cz^k,\cz^k).$$
Any such $B$ is called a parametrix of $A$. 
\end{thm}

Note that parametrices of elliptic operators are uniquely determined modulo smoothing operators. 
Recall once more that smoothing operators are precisely those integral operators with an integral kernel 
which is rapidly decreasing in all variables. 

%%%%%%%%%%%%%%%%%%%%%%%%%%%%%%%%%%%%%%%%%%%%%%%%%%
\subsection{Ellipticity and Fredholm property}\label{sec:02.3}

Let $A\in \Gamma^{\mu,\nu}(\rz^{m}\times \rz^{n};\cz^{k},\cz^{k})$. If $A$ is elliptic one can construct a 
parametrix $B\in\Gamma^{-\mu,-\nu}(\rz^{m}\times \rz^{n};\cz^{k},\cz^{k})$, i.e.,  
both $1-AB$ and $1-BA$ are smoothing operators. Since smoothing operators induce compact 
operators in the Sobolev spaces of any order, the implication a$)\Rightarrow$ b$)$ of 
the following theorem is evident:

\begin{thm}\label{thm:shubin}
For $A\in \Gamma^{\mu,\nu}(\rz^{m}\times \rz^{n};\cz^{k},\cz^{k})$ the following properties are equivalent$:$
\begin{itemize}
 \item[a$)$] $A$ is elliptic.  
 \item[b$)$] For every $(s,t)\in\rz^2$, $A$ induces Fredholm operators 
   $$Q^{s,t}(\rz^{m}\times \rz^{n};\cz^{k})\lra 
       Q^{s-\mu,t-\nu}(\rz^{m}\times \rz^{n};\cz^{k}).$$   
 \item[c$)$] There exists a tuple $(s,t)\in\rz^2$ such that $A$ induces a Fredholm operator 
   $$Q^{s,t}(\rz^{m}\times \rz^{n};\cz^{k})\lra 
       Q^{s-\mu,t-\nu}(\rz^{m}\times \rz^{n};\cz^{k}).$$   
\end{itemize}
\end{thm}

{
The implication b$)\Rightarrow$ c$)$ is trivial. 
In the sequel we shall prove the implication c$)\Rightarrow$ a$)$. 
%To shorten notation let us now assume $k=1$. By using order reductions we also may assume without 
%loss of generality that $\mu=\nu=s=t=0$, i.e., we may assume that 
%$A\in \Gamma^{0,0}(\rz^{m}\times \rz^{n})$ is a Fredholm operator in $L^2(\rz^{m}\times \rz^{n})$. 
The method of proof is inspired by that of Theorem 1 in Section 
2.3.4.1 of \cite{ReSc} and by that of Theorem 1.6 in \cite{ScSe}. 
}

%%%%%%%%%%%%%%%%%%%%%%%%%%%%%%%%%%%%%%%%%%%%%%%%%%
\subsubsection{A family of isometries}\label{sec:02.3.1}

Let $E$ be a Hilbert space. For fixed $(x_0,\xi_0)\in\rz^{n}\times\rz^n$ with $|(x_0,\xi_0)|=1$ and 
{an arbitrarily fixed $\tau\in(0,1/2)$} define $S_\lambda\in \scrL(L^2(\rz^n,E))$, $\lambda\ge1$, by 
\begin{equation}\label{eq:slambda}
 (S_\lambda u)(x)=\lambda^{n\tau/2}e^{i\lambda x\xi_0}u\big(\lambda^\tau(x-\lambda x_0)\big).
\end{equation}
It is straightforward to verify that any $S_\lambda$ is an isometric isomorphism with inverse given by 
 $$(S_\lambda^{-1} v)(x)=\lambda^{-n\tau/2}e^{-i\lambda(\lambda x_0+\lambda^{-\tau}x)\xi_0}
     v\big(\lambda^\tau(\lambda x_0+\lambda^{-\tau}x)\big).$$
Moreover, 
\begin{align}\label{eq:weak}
 \mathop{\mbox{\textrm{w-}$\lim$}}_{\lambda\to+\infty} S_\lambda u =0\qquad \forall\;u\in L^2(\rz^n,E), 
\end{align}
where $\mathop{\mbox{\textrm{w-}$\lim$}}$ denotes the limit with respect to the weak topology of 
$L^2(\rz^n,E)$. In fact, this property follows from the fact that all $S_\lambda$ are isometries and that 
\begin{align*}
 |(S_\lambda u,v)_{L^2(\rz^n,E)}|&=\Big|\int \big(S_\lambda u(x),v(x)\big)_E\,dx\Big|\\
 &\le \int \lambda^{n\tau/2}\|u(\lambda^\tau(x-\lambda x_0))\|_E\|v(x)\|_E\,dx\\
 &\le \lambda^{-n\tau/2}\|u\|_{L^1(\rz^n,E)}\|v\|_{L^\infty(\rz^n,E)}\xrightarrow{\lambda\to+\infty}0
\end{align*} 
for every $u$ and $v$ belonging to the dense subspace $\scrS(\rz^n,E)$ of $L^2(\rz^n,E)$. 

%%%%%%%%%%%%%%%%%%%%%%%%%%%%%%%%%%%%%%%%%%%%%%%%%%
\subsubsection{Recovering the principal symbol}\label{sec:02.3.2}

Let $a\in \Gamma^{\nu}(\rz^n,\scrL(E))$ be an operator-valued symbol in the sense of Section \ref{sec:02.1.1} 
{
For convenience of notation we assume that $a$ is $\scrL(E)$-valued, but the following results remain 
valid for the more general case of $a$ being $\scrL(E,F)$-valued, with two Hilbert spaces $E$ and $F$.
}
If the $S_\lambda$, $\lambda\ge 1$, are as introduced in the previous Section \ref{sec:02.3.1}, 
a direct calculation shows that 
\begin{align}\label{eq:conjugation}
 S_\lambda^{-1}\op(a) S_\lambda=\op(a_\lambda),\qquad 
 a_\lambda(x,\xi)=a(\lambda x_0+\lambda^{-\tau}x,\lambda\xi_0+\lambda^\tau\xi). 
\end{align}
Note that $a_\lambda\in \Gamma^{\nu}(\rz^n,\scrL(E))$ for every $\lambda$. The following estimate will be crucial 
later on$:$

\begin{lem}\label{lem:key-estimate}
Let $a\in \Gamma^{\nu}(\rz^n,\scrL(E))$ with ${\nu}\le 0$ and $\rho=\frac{\tau}{1-\tau}$ $($note that $0<\rho<1)$. 
Then, for any order of derivatives, 
  $$\big \|D^\alpha_\xi D^\beta_x a_\lambda(x,\xi) \big\|_{\scrL(E)}
      \le C_{\alpha\beta}\, \lambda^{(1-\tau){\nu}-\tau|\beta|}\,\spk{x,\xi}^{\rho|\alpha|-{\nu}}$$
uniformly in $(x,\xi)\in\rz^{n}\times\rz^n$ and $\lambda\ge1$. 
\end{lem}
\begin{proof}
By chain rule and using the standard symbol estimates for $a$, we have 
  $$\big \|D^\alpha_\xi D^\beta_x a_\lambda(x,\xi) \big\|_{\scrL(E)}
      \le C\,\lambda^{|\alpha|\tau-|\beta|\tau}\,
      \spk{\lambda x_0+\lambda^{-\tau}x,\lambda\xi_0+\lambda^\tau\xi}^{{\nu}-\rho|\alpha|},$$
with a constant $C$ independent of $(x,\xi)$ and $\lambda$. 
Since $\spk{v+w}^{-1}\le C\spk{w}/|v|$ by Peetre's inequality and $\spk{\sigma w}\le\sigma\spk{w}$ 
for $\sigma\ge 1$, we can estimate 
\begin{align*}
 \spk{\lambda x_0+\lambda^{-\tau}x,\lambda\xi_0+\lambda^\tau\xi}^{{\nu}-\rho|\alpha|}
 &\le C\spk{\lambda^{-\tau} x,\lambda^{\tau}\xi}^{\rho|\alpha|-{\nu}}
    |(\lambda x_0,\lambda\xi_0)|^{{\nu}-\rho|\alpha|}\\
 &\le C\lambda^{(\rho|\alpha|-{\nu})\tau}\lambda^{{\nu}-\rho|\alpha|}\spk{x,\xi}^{\rho|\alpha|-{\nu}}, 
\end{align*}
resulting in 
  $$\big \|D^\alpha_\xi D^\beta_x a_\lambda(x,\xi) \big\|_{\scrL(E)}
      \le C\,\lambda^{(1-\tau){\nu}-\tau|\beta|+(\tau-\rho+\tau\rho)|\alpha|}\spk{x,\xi}^{\rho|\alpha|-{\nu}}.$$
It remains to observe that $\tau-\rho+\tau\rho=0$, due to the choice of $\rho$. 
\end{proof}

\begin{lem}\label{lem:lebesgue}
Let $\{a_\lambda\mid\lambda\ge1\})$ be a subset of $\Gamma^0(\rz^n,\scrL(E))$, 
$\sigma\in\cz$ a constant, and $u\in\scrS(\rz^n,E)$. Assume that 
\begin{itemize}
 \item[(1a)] $a_\lambda(x,\xi)\xrightarrow{\lambda\to+\infty}\sigma$ for all $(x,\xi)\in\rz^n\times\rz^n$, 
 \item[(1b)] for every $x\in\rz^n$ there exist constants $c_x,m_x\ge 0$ such that 
   $$\|a_\lambda(x,\xi)\|\le c_x\spk{\xi}^{m_x}\qquad \forall\;\xi\in\rz^n\quad\forall\;\lambda\ge1,$$ 
 \item[(2)] there exists a $g\in L^1(\rz^n)$ such that 
  $$\|[\op(a_\lambda)u](x)\|_E^2\le g(x)\qquad \forall\; x\in\rz^n\quad\forall\;\lambda\ge 1.$$
\end{itemize}
Then $\op(a_\lambda)u\xrightarrow{\lambda\to+\infty}\sigma u$ in $L^2(\rz^n,E)$. 
\end{lem}
\begin{proof}
The result follows directly from Lebegue's dominated convergence theorem, provided we can show that 
$\op(a_\lambda)u$ converges pointwise on $\rz^n$ to $\sigma u$ as $\lambda$ tends to infinity. 
However, with $x\in\rz^n$ fixed, 
 $$[\op(a_\lambda)u](x)=\int e^{ix\xi}a_\lambda(x,\xi)\wh{u}(\xi)\,\dbar\xi.$$
By assumption (1a), the integrand converges pointwise on $\rz^n_\xi$ to $\sigma e^{ix\xi}\wh{u}(\xi)$. 
By (1b) the integrand is majorized in norm by $h(\xi):=c_x\spk{\xi}^{m_x}\wh{u}(\xi)\in L^1(\rz^n_\xi)$.  
Thus, by dominated convergence, 
 $$[\op(a_\lambda)u](x)\xrightarrow{\lambda\to+\infty}
     \sigma\int e^{ix\xi}\wh{u}(\xi)\,\dbar\xi=\sigma u(x).$$
This completes the proof. 
\end{proof}

The following proposition gives a method for recovering the principal symbol from the operator$:$

\begin{prop}\label{prop:convergence}
Let $A=\op(a)\in \Gamma^0_\cl(\rz^n,\scrL(E))$, $a_\lambda$ as in \eqref{eq:conjugation}, and $u\in\scrS(\rz^n,E)$. 
Then 
 $$\op(a_\lambda)u\xrightarrow{\lambda\to+\infty}a^{(0)}(x_0,\xi_0)u\quad\text{in }L^2(\rz^n,E),$$
where $a^{(0)}\in\scrC^\infty(\mathbb{S}^{2n-1},\scrL(E))$ denotes the homogeneous principal symbol of $a$. 
\end{prop}
\begin{proof}
By Lemma \ref{lem:key-estimate} with $|\alpha|=|\beta|={\nu}=0$,  condition (1b) of Lemma \ref{lem:lebesgue} 
is obviously satisfied $($with $m_x=0)$. 
Now let $\chi(x,\xi)$ be a zero-excision function and write $a=a^0+r$, where 
 $$a^0(x,\xi)=\chi(x,\xi) a^{(0)}(x,\xi),\qquad r\in \Gamma^{-1}(\rz^n,\scrL(E)).$$
Then $a_\lambda=a^0_\lambda+r_\lambda$. By Lemma \ref{lem:key-estimate} with $|\alpha|=|\beta|=0$ and 
$\mu=-1$, it is clear that $r_\lambda(x,\xi)\to0$ for all $x$ and $\xi$. Moreover, by homogeneity of $a^{(0)}$, 
\begin{align*}
 a^0_\lambda(x,\xi)
 =&\chi(\lambda x_0+\lambda^{-\tau}x,\lambda\xi_0+\lambda^\tau\xi)
      a^{(0)}(x_0+\lambda^{-1-\tau}x,\xi_0+\lambda^{-1+\tau}\xi)
\end{align*}
and thus $a^0_\lambda(x,\xi)\to a^{(0)}(x_0,\xi_0)$ for all $x$ and $\xi$. Therefore assumption (1a) 
of Lemma \ref{lem:lebesgue} with $\sigma=a^{(0)}(x_0,\xi_0)$ is satisfied. 

It remains to verify assumption (2). To this end let $M\in\nz$ and write, using integration by parts,  
 $$\spk{x}^{2M}[\op(a_\lambda)u](x)
     =\int e^{ix\xi}(1+\Delta_\xi)^M\big(a_\lambda(x,\xi)\wh{u}(\xi)\big)\,\dbar\xi.$$
By product rule and Lemma \ref{lem:key-estimate} there exist functions $u_\alpha\in\scrS(\rz^n,E)$ such that 
 $$\spk{x}^{2M}\|[\op(a_\lambda)u](x)\|_E\le
     \sum_{|\alpha|\le 2M}\int \spk{x,\xi}^{\rho|\alpha|}\wh{u_\alpha}(\xi)\,\dbar\xi.$$
Hence 
 $$\|[\op(a_\lambda)u](x)\|^2_E\le C\spk{x}^{4M(\rho-1)}=:g(x)$$
with a suitable constant independent of $x$ and $\lambda$. 
Since $\rho-1<0$ we can choose $M$ so large that $g\in L^1(\rz^n)$.  
\end{proof}

%%%%%%%%%%%%%%%%%%%%%%%%%%%%%%%%%%%%%%%%%%%%%%%%%%
\subsubsection{The proof of Theorem \textnormal{\ref{thm:shubin}}}\label{sec:02.3.3}

First we shall proof the following result on pseudodifferential operators with operator-valued 
symbols. {Recall that a linear continuous operator is called upper semi-fredholm if it has closed range and 
finite-dimensional kernel; it is called lower semi-fredholm, if its range is closed and of finite co-dimension$:$}

\begin{prop}\label{prop:semi-fred}
Consider $A=\op(a)\in \Gamma^0_\cl(\rz^n,\scrL(E))$ as a bounded operator in $L^2(\rz^n,E)$ and let 
$(x_0,\xi_0)\in\rz^n\times\rz^n$ be a unit-vector.  
\begin{itemize}
 \item[a$)$] If $A$ is upper semi-fredholm, $a^{(0)}(x_0,\xi_0)$ is injective. 
 \item[b$)$] If $A$ is lower semi-fredholm, $a^{(0)}(x_0,\xi_0)$ is surjective.
\end{itemize}
\end{prop} 
\begin{proof}
Assume that $A=\op(a)\in \Gamma^0_\cl(\rz^n,\scrL(E))$ induces an upper semi-fredholm operator 
$A\in \scrL(L^2(\rz^n,E))$. Since $E$ is a Hilbert space, there exists a $B\in\scrL(L^2(\rz^n,E))$ 
such that $K:=1-BA$ is a compact operator in $L^2(\rz^n,E)$. 

Let $u\in\scrS(\rz^n)$ with $\|u\|_{L^2(\rz^n)}=1$ and define 
$u_e\in\scrS(\rz^n,E)$, $e\in E$, by $u_e(x)=u(x)e$. Then, with notations from the previous subsection, 
\begin{align*}
 \|e\|_E
 =&\|u_e\|_{L^2(\rz^n,E)}=\|(BA+K)S_\lambda u_e\|_{L^2(\rz^n,E)}\\
 \le& \|B\|_{\scrL(L^2(\rz^n,E))}\|S_\lambda^{-1}AS_\lambda u_e\|_{L^2(\rz^n,E)}+
          \|K S_\lambda u_e\|_{L^2(\rz^n,E)}\\
     &\xrightarrow{\lambda\to+\infty}\|B\|_{\scrL(L^2(\rz^n,E))}\|a^{(0)}(x_0,\xi_0)e\|_E.
\end{align*}
For the convergence we have used that $KS_\lambda u_e\to 0$, since $S_\lambda u_e\to 0$ weakly by 
\eqref{eq:weak} and $K$ is compact, and that $S_\lambda^{-1}AS_\lambda u_e=\op(a_\lambda)u_e\to u_e$ in 
$L^2(\rz^n,E)$ due to Proposition \ref{prop:convergence}. 
Therefore, 
 $$\|a^{(0)}(x_0,\xi_0)e\|_E\ge \frac{1}{\|B\|_{\scrL(L^2(\rz^n,E))}}\|e\|_E\qquad\forall\;e\in E.$$
This implies a$)$. 
If $A$ is lower semi-fredholm, its adjoint is an upper semi-fredholm operator. By a$)$, the principal symbol of $A^*$ 
evaluated in $(x_0,\xi_0)$, i.e., $a^{(0)}(x_0,\xi_0)^*$, is injective. Hence $a^{(0)}(x_0,\xi_0)$ is surjective. 
\end{proof}

{
Let us emphasize once more that the previous result remains valid in case of 
$A=\op(a)\in\Gamma^0_\cl(\rz^n,\scrL(E,F))$ with Hilbert spaces $E$ and $F$, considered as an operator from 
$L^2(\rz^n,E)$ to $L^2(\rz^n,F)$. 
}

The proof of c$)\Rightarrow$ a$)$ of Theorem \ref{thm:shubin} now 
{works as follows:
Consider $A\in \Gamma^{\mu,\nu}_\cl(\rz^{m}\times \rz^{n};\cz^{k\times k})$ as an operator with operator-valued 
symbol $a\in \Gamma^\nu_\cl(\rz^n,\scrL(E,F))$ with $E=Q^s(\rz^m,\cz^k)$ and $F=Q^{s-\mu}(\rz^m,\cz^k)$. 
With order-reductions $\Lambda^{s}_E=\op(\lambda^s_E)$  and $\Lambda^{s}_F=\op(\lambda^s_F)$ as 
described in Section \ref{sec:02.1.3}, using \eqref{eq:redord} define $\wt{A}:=\Lambda^{t-\nu}_F A \Lambda^{-t}_E$. Then 
$\wt{A}=\op(\wt{a})\in\Gamma^{0}_\cl(\rz^n,\scrL(E,F))$ and the Fredholm property of $A$ is equivalent to that of 
$\wt{A}:L^2(\rz^n,E)\to L^2(\rz^n,F)$. By Proposition \ref{prop:semi-fred}, the homogeneous principal symbol 
$\wt{a}^{(0)}\in\scrC^\infty(\mathbb{S}^{2n-1},\scrL(E,F))$ is pointwise invertible. However, this principal symbol 
just coincides with $\sigma^\nu_2(A)$ as introduced in 
Section \ref{sec:02.2.2}. Analogously, $\sigma^\mu_1(A)$ evaluated in an arbitrary unit-vector of $\rz^m\times\rz^m$ 
is invertible as an operator in $Q^t(\rz^n,\cz^k)\to Q^{t-\nu}(\rz^n,\cz^k)$. 
}

{
\begin{rem}
Let us mention an alternative approach to prove Theorem \ref{thm:shubin}, based on $C^*$-algebraic arguments. 
Let $\Gamma(\rz^n)$ denote the $C^*$-closure of $\Gamma^0_\cl(\rz^n)$ and $\mathcal{K}_n$ the space of 
compact operators in $L^2(\rz^n)$. Then $\Gamma(\rz^n)/\mathcal{K}_n$ can be identified with the space of 
continuous functions  on the unit-sphere $\mathbb{S}^{2n-1}$; see \cite{Bohl} for details. Using \eqref{eq:tensor}, the 
$C^*$-closure of $\Gamma^{0,0}(\rz^m\times\rz^n)$, factored by the compact operators, can be identified with 
$\big[(\Gamma(\rz^m)/\mathcal{K}_m)\otimes\Gamma(\rz^n)\big]\oplus
\big[\Gamma(\rz^m)\otimes (\Gamma(\rz^n)/\mathcal{K}_n)\big]$. This means that an operator 
$($from the $C^*$-closure$)$ is Fredholm if, and only if, the two associated principal symbols are invertible. 
Filling in the details of the above argument is of a complexity comparable with that of the proof above . 
\end{rem}
}

%%%%%%%%%%%%%%%%%%%%%%%%%%%%%%%%%%%%%%%%%%%%%%%%%%
\subsection{Spectral invariance}\label{sec:02.4} 

A consequence of Theorem \ref{thm:shubin} is the following result, the so-called \emph{spectral-invariance} 
of bisingular pseudodifferential operators$:$  

\begin{thm}\label{thm:spectral}
Let $A\in \Gamma^{\mu,\nu}(\rz^{m}\times \rz^{n};\cz^{k},\cz^{k})$. Assume that $A$ induces an isomorphism 
$Q^{s,t}(\rz^{m}\times \rz^{n};\cz^{k})\lra Q^{s-\mu,t-\nu}(\rz^{m}\times \rz^{n};\cz^{k})$ for some tuple 
$(s,t)\in\rz^2$. Then there exists a $B\in \Gamma^{\mu,\nu}(\rz^{m}\times \rz^{n};\cz^{k},\cz^{k})$ such that 
$AB=BA=1$. 
In particular, $A$ induces an isomorphism 
$Q^{s,t}(\rz^{m}\times \rz^{n};\cz^{k})\lra Q^{s-\mu,t-\nu}(\rz^{m}\times \rz^{n};\cz^{k})$ for every 
tuple $(s,t)\in\rz^2$.
\end{thm}

In other words, invertibility as a bounded operator between Sobolev spaces implies the invertibility within the 
class of bisingular pseudodifferential operators. 

\begin{proof}
To shorten notation let us assume $k=1$. The isomorphism is, in particular, a Fredholm operator. 
Due to Theorem \ref{thm:shubin}, $A$ is elliptic. Therefore 
it has a parametrix $B_0\in \Gamma^{-\mu,-\nu}(\rz^{m}\times \rz^{n})$. Thus $K_R:=1-AB_0$ 
and $K_L:=1-B_0A$ are smoothing operators. Passing to the action in Sobolev spaces, and resolving both equations 
for $A^{-1}$ we 
{obtain $A^{-1}=A^{-1}K_R+B_0$ and $A^{-1}=K_LA^{-1}+B_0$. Inserting the latter equation 
in the previous one yields} 
 $$A^{-1}=B_0+B_0K_R+K_LA^{-1}K_R.$$
Obviously, both $B_0$ and $B_0K_R$ belong to $\Gamma^{-\mu,-\nu}(\rz^{m}\times \rz^{n})$. 
Now let $R:=K_LA^{-1}K_R$. We shall argue below that $R$ is smoothing and 
therefore $B=B_0+B_0K_R+R\in \Gamma^{-\mu,-\nu}(\rz^{m}\times \rz^{n})$ is the desired operator.

Since $K_L$ and $K_R$ are smoothing it is obvious that both $R$ and $R^*$ map 
$L^2(\rz^{m}\times \rz^{n})$ to $\scrS(\rz^{m}\times \rz^{n})$. However, this is known to be equivalent 
to $R$ being an integral operator with an integral kernel that is rapidly decreasing in all variables; for convenience 
of the reader we sketch the argument: First of all one sees that $R$ has a kernel 
$k(x,y)=k(x_1,x_2,y_1,y_2)\in L^2(\rz^{2m}_x\times\rz^{2n}_y)$ such that 
 $$k\in\scrS(\rz^n_{x_1}\times \rz^{m}_ {y_1},L^2(\rz^n_{x_2}\times \rz^{m}_ {y_2}))
     \cap 
     \scrS(\rz^n_{x_2}\times \rz^{m}_ {y_2},L^2(\rz^n_{x_1}\times \rz^{m}_ {y_1})).$$
Thus the claim follows if we can show that 
 $$\scrS(\rz^k_u,L^2(\rz^\ell_v))\cap \scrS(\rz^\ell_v,L^2(\rz^k_u))=\scrS(\rz^{k+\ell}_{(u,v)}).$$
Let $g$ be a function from the space on the left-hand side and denote by $\|\cdot\|$ the norm of $L^2(\rz^{k+\ell})$. 
Then, by Parseval's identity, 
 $$\|g\|=(2\pi)^{-(k+\ell)/2}\|\scrF g\|=(2\pi)^{-k/2}\|\scrF_{u\to\xi}g\|=(2\pi)^{-\ell/2}\|\scrF_{v\to\eta}g\|.$$
Combining this repeatedly with the estimate $ab\le a^2+b^2$, one obtains that 
\begin{align*}
 \|\spk{u}^{i}\spk{v}^{j}\spk{D_u}^{i^\prime}\spk{D_v}^{j^\prime}g\|
 \le & \;
    C\Big(\|\spk{u}^{4i}\spk{D_u}^{i^\prime}g\|+\|\spk{D_u}^{2i^\prime}g\|+\|\spk{D_u}^{4i^\prime}g\|+\\
    &+\|\spk{v}^{4j}\spk{D_v}^{j^\prime}g\|+\|\spk{D_v}^{2j^\prime}g\|+\|\spk{D_v}^{4j^\prime}g\|\Big)
\end{align*}
is finite for any choice of non negative integers $i,i^\prime,j,j^\prime$. 
This yields that $g$ belongs to $\scrS(\rz^{k+\ell})$. 
\end{proof}

\begin{cor}\label{cor:extension}
Let $A\in \Gamma^{\mu,\nu}(\rz^{m}\times \rz^{n};\cz^{k},\cz^{k})$ be elliptic and $\mu,\nu\ge 0$. 
Then the unbounded operator 
 $$A_{s,t}:\scrS(\rz^{m}\times \rz^{n},\cz^k)\subset 
     Q^{s,t}(\rz^{m}\times \rz^{n},\cz^k)\lra Q^{s,t}(\rz^{m}\times \rz^{n},\cz^k)$$
has one, and only one, closed extension, given by the action of $A$ on the domain 
$Q^{s+\mu,t+\mu}(\rz^{m}\times \rz^{n},\cz^k)$. 
The spectrum of the closure of $A_{s,t}$ does not depend on both $s$ and $t$.  
\end{cor}
\begin{proof}
By density of the rapidly decreasing functions in any Sobolev space, it is clear that 
$Q^{s+\mu,t+\mu}(\rz^{m}\times \rz^{n},\cz^k)$ is contained in the domain of the closure of $A_{s,t}$. 
Moreover, if both $u$ and $Au$ belong to $Q^{s,t}(\rz^{m}\times \rz^{n},\cz^k)$ then 
$u\in Q^{s+\mu,t+\mu}(\rz^{m}\times \rz^{n},\cz^k)$ by elliptic regularity. Therefore, the domain of 
any closed extension is a subset of, and hence equal to, $Q^{s+\mu,t+\mu}(\rz^{m}\times \rz^{n},\cz^k)$. 

The statement on the spectrum follows directly from Theorem \ref{thm:spectral} and the fact that 
$\lambda-A\in \Gamma^{\mu,\nu}(\rz^{m}\times \rz^{n};\cz^{k},\cz^{k})$ for any $\lambda\in\cz$.
\end{proof}

%%%%%%%%%%%%%%%%%%%%%%%%%%%%%%%%%%%%%%%%%%%%%%%%%%
\section{Bisingular operators on closed manifolds}\label{sec:03} 

In \cite{Rodi} bisingular operators acting on sections in vector bundles over products of closed manifolds are 
considered. We shall use the notation $L^{\mu,\nu}_{\cl}(M\times N; E,F)$ for such operators and 
$Q^{s,t}(M\times N,G)$ for the associated Sobolev spaces, where $M$ and $N$ are closed Riemannian manifolds 
and $E$, $F$ and $G$ are finite-dimensional hermitian vector-bundles over $M\times N$.
%\footnote{A comment analogous to the one of footnote \ref{fn:01} applies also here.}

%%%%%%%%%%%%%%%%%%%%%%%%%%%%%%%%%%%%%%%%%%%%%%%%%%
\subsection{Description of the calculus}\label{sec:03.1} 

As usual, bisingular operators on a manifold are defined as those that in any local trivialisation of the bundles and any 
local coordinates correspond to bisingular operators in a product of two Euclidean spaces, with symbols taking values 
in $\cz^{\mathrm{dim}\,F\times\mathrm{dim}\,E}$. We shall not go too much into the details, but only describe 
how the classes $\Gamma^{\mu,\nu}$ introduced above have to be modified to recover the situation of \cite{Rodi}. 

%%%%%%%%%%%%%%%%%%%%%%%%%%%%%%%%%%%%%%%%%%%%%%%%%%
\subsubsection{The calculus on $\rz^m\times\rz^n$}\label{sec:03.1.1} 

For a Fr\'{e}chet space $F$ define the space $L^\nu(\rz^n,F)$ as in the beginning of Section \ref{sec:02.1}, 
replacing in \eqref{eq:symbol1} the term $\spk{x,\xi}^{|\alpha|+|\beta|-\nu}$ by $\spk{\xi}^{|\alpha|-\nu}$. 

For defining the classical symbols $L^\nu_\cl(\rz^n,F)$, in the subsequent part one considers homogeneous 
components $a^{(\nu-j)}:\rz^n\times(\rz^n_\xi\setminus\{0\})\to F$ which are homogeneous in the sense of 
 $$a^{(\nu-j)}(x,t\xi)=t^{\nu-j}\,a^{(\nu-j)}(x,\xi)\qquad \forall\;t>0\quad\forall\;x\quad\forall\;\xi\not=0.$$
The excision function $\chi(x,\xi)$ needs to be replaced by an excision function $\chi(\xi)$. 

Starting out with these symbol classes, one then introduces, as before, the bisingular symbols 
$L^{\mu,\nu}_{\cl}(\rz^m\times\rz^n;\cz^k,\cz^\ell)$. The corresponding Sobolev spaces 
$Q^{s,t}(\rz^m\times\rz^n)$ are defined as the closure of $\scrS(\rz^m\times\rz^n)$ with respect to the norm 
$\|u\|_{s,t}=\|\Lambda^{s,t}u\|_{L^2(\rz^m\times\rz^n)}$, where  $\Lambda^{s,t}$ is the operator with 
symbol $\lambda^{s,t}(\xi,\eta)=\spk{\xi}^s\spk{\eta}^t$. 

The two principal symbols associated with 
$A=\op(a)\in L^{\mu,\nu}_\cl(\rz^m\times\rz^n;\cz^k,\cz^\ell)$ are then  
\begin{align}\label{eq:symb}
\begin{split}
  \sigma_1^\mu(A)&=a^{(\mu)}_1\in 
    \scrC^\infty\big(\rz^m_x\times \mathbb{S}^{m-1}_\xi,L^\nu_\cl(\rz^n;\cz^{\ell\times k})\big), \\ 
  \sigma_2^\nu(A)&=a^{(\nu)}_2\in 
    \scrC^\infty\big(\rz^n_y\times \mathbb{S}^{n-1}_\eta,L^\mu_\cl(\rz^m;\cz^{\ell\times k})\big),   
\end{split}
\end{align}
and ellipticity asks the pointwise invertibility of both these symbols. 

The analogue of Theorem \ref{thm:parametrix} holds true, while Theorem \ref{thm:shubin} fails to be true, 
since smoothing operators do not induce compact operators in the Sobolev spaces of $\rz^m\times\rz^n$. 
However, the analogue of Theorem \ref{thm:shubin} for operators on a product of compact manifolds is valid, 
as we shall see below. 

%%%%%%%%%%%%%%%%%%%%%%%%%%%%%%%%%%%%%%%%%%%%%%%%%%
\subsubsection{The principal symbols}\label{sec:03.1.2} 

For an operator $A\in L^{\mu,\nu}_{\cl}(M\times N; E,F)$ the existence of local principal symbols leads to 
two globally defined $($on the unit co-sphere bundles $S^*M$ and $S^*N$, respectively$)$ objects, 
again denoted by $\sigma_1^\mu(A)$ and $\sigma_2^\nu(A)$. 
If $v=(x,\xi)\in S^*M$ then $\sigma^\mu_1(A)(v)$ is an operator in 
$L^\nu_\cl(N;E({x}),F({x}))$, where $L^\nu_\cl$ refers to the usual space of classical pseudodifferential 
operators on a closed manifold and 
 $$E(x):=E\big|_{\{x\}\times N},\quad F(x):=F\big|_{\{x\}\times N}\qquad x\in M,$$
considered as vector bundles over $N\cong\{x\}\times N$. 

If we denote by $\pi_M:S^*M\to M$ the canonical projection and define the (infinite-dimensional) Hilbert space 
bundle $\calQ^s(N,E)$ over $M$ by taking as fibre in $m\in M$ the Sobolev space 
$Q^s(N,E(m))$ of sections in $E(m)$ $($see Section \ref{sec:05} for details$)$\footnote{{The common 
notation for these Sobolev spaces is $H^s$; however, for reasons of consistency with the previously 
employed notation we shall use the letter $Q$ rather than $H$.}}, 
then we can consider $\sigma^\mu_1(A)$ as a bundle homomorphism 
\begin{align}\label{eq:hom1}
 \sigma^\mu_1(A):\pi_M^*\calQ^s(N,E)\lra \pi_M^*\calQ^{s-\nu}(N,F),\qquad s\in\rz.
\end{align}
Similarly, 
\begin{align}\label{eq:hom2}
 \sigma^\nu_2(A):\pi_N^*\calQ^s(M,E)\lra \pi_N^*\calQ^{s-\mu}(M,F),\qquad s\in\rz.
\end{align}

\begin{thm}
$A\in L^{\mu,\nu}_{\cl}(M\times N; E,F)$ is called elliptic if both homomorphisms \eqref{eq:hom1} and 
\eqref{eq:hom2} are isomorphisms\footnote{{Evaluation of the principal symbols in a specific co-vector 
gives a standard, classical pseudodifferential operator of order $\mu$ respectively $\nu$ on the manifold $M$ or $N$, 
respectively.} Due to spectral invariance of this calculus, {conditions \eqref{eq:hom1} and \eqref{eq:hom2}} 
are independent of $s$.}. 
Then, the following are equivalent$:$
\begin{itemize}
 \item[a$)$] $A\in L^{\mu,\nu}_{\cl}(M\times N; E,F)$ is elliptic. 
 \item[b$)$] There exists a $B\in L^{-\mu,-\nu}_{\cl}(M\times N; F,E)$ such that both $1-AB$ and $1-BA$ are 
  smoothing operators. 
\end{itemize}
\end{thm}

%%%%%%%%%%%%%%%%%%%%%%%%%%%%%%%%%%%%%%%%%%%%%%%%%%
\subsection{Ellipticity and Fredholm property}\label{sec:03.2} 

We are now going to explain that the analogue of Theorem \ref{thm:shubin} holds for operators 
$A\in L^{\mu,\nu}_{\cl}(M\times N; E,F)$. 
{
Assume that $A$ induces a Fredholm operator 
 $$A:Q^{s,t}(M\times N,E)\lra Q^{s-\mu,t-\nu}(M\times N,F)$$
for some fixed numbers $s$ and $t$. Let $B$ be the corresponding inverse modulo compact operators. 
}
Let  $K:=1-BA$ and $v_0\in S^*_{m_0}M$ be a given, fixed unit co-vector. We shall verify the invertibility of 
 $${\sigma^{\mu}_1(A)(v_0)\in L^\nu_\cl(N;E(m_0),F(m_0)).}$$
To this end, let  $U$ be a coordinate system of $M$ near 
$m_0$ such that $v_0$ corresponds to $(x_0,\xi_0)$ and that $E|_{U\times N}\cong U\times E(m_0)$, 
$F|_{U\times N}\cong U\times F(m_0)$ in the sense of Proposition \ref{prop:triv}. 
Moreover, let $\chi_1,\chi_2,\chi_3\in\scrC^\infty_0(U_0)$ such that 
$\chi_{i+1}\equiv 1$ on the support of $\chi_i$ for $i=1,2$. Consider the $\chi_i$ as functions on $M\times N$, 
not depending on the variable of $N$. Multiplying the identity $K=1-BA$ from the left with $\chi_1$ , from the 
right with $\chi_3$, and rearranging terms yields 
\begin{equation}\label{eq:local}
 \chi_1 B\chi_2\chi_3 A\chi_3=\chi_1-\chi_1K\chi_3-\chi_1B(1-\chi_2)A\chi_3.
\end{equation}
Note that $(1-\chi_2)A\chi_3\in L^{-\infty,{\nu}}_{\cl}(M\times N; E,F)$ due to the disjoint supports of 
$(1-\chi_2)$ and $\chi_3$, {and that all four operators in \eqref{eq:local}  
%$B^\prime:=\chi_1 B\chi_2$, $A^\prime:=\chi_3 A\chi_3$, $K_1^\prime:=\chi_1K\chi_3$ and 
%$K_2^\prime:=\chi_1B(1-\chi_2)A\chi_3$ 
are localized in $U\times N$. 
In particular, they can be identified -- after passing to local coordinates in $U$ -- with operators on 
$\rz^m\times N$. 

Now let $\lambda^s(\xi)=[\xi]^s$, $s\in\rz$, where $[\cdot]$ denotes a smooth, positive function that coincides with 
the usual modulus outside some neighborhood of the origin. Obviously, $\lambda^s\in L^s_\cl(\rz^m)$ and 
$\op(\lambda^s)\op(\lambda^{-s})=1$. 
Define the operators $\Lambda^s=\mathrm{op}(\lambda^s)\otimes 1$, $s\in\rz$,  on $\rz^m\times N$. 

Multiplying \eqref{eq:local} from the left with $\Lambda^s$, from the right with $\Lambda^{-s}$, and by substituting 
on the left-hand side $\chi_2\chi_3$ by $\chi_2\Lambda^{\mu-s}\Lambda^{s-\mu}\chi_3$, we obtain an equality 
 $$B^\prime A^\prime=\Phi-K_1-K_2,$$ 
with obvious meaning of notation. In particular, $A^\prime$ and $\Phi$ are pseudodifferential operators with 
respective operator-valued symbols 
\begin{align*}
 a & \in L^{0}_\cl(\rz^m,\scrL({Q^t}(N,E(m_0)),{Q^{t-\nu}}(N,F(m_0))),\\
 \varphi & \in L^{0}_\cl(\rz^m,\scrL({Q^t}(N,E(m_0)),{Q^{t}}(N,F(m_0))),
\end{align*}  
where $a^{(0)}(x_0,\xi_0)$ is the local expression of $\sigma^0_1(A)(v_0)$ and $\varphi^{(0)}(x_0,\xi_0)=1$. 

Observe that $K_{{2}}$ is not a compact operator, but extends to a continuous map 
$L^1(\rz^m,Q^t(N,E(m_0)))$ into $L^2(\rz^m,Q^{t}(N,F(m_0)))$. 
The injectivity of $\sigma^0_1(A)(v_0)$ now follows from the following proposition; its 
surjectivity, hence invertibility, then follows by considering the adjoint of $A$.   

\begin{prop}
Let $E,F$ be two Hilbert spaces and $(x_0,\xi_0)\in \rz^m\times\rz^m$ with $|\xi_0|=1$.  
Moreover let $A=\op(a)\in L^0_\cl(\rz^m,\scrL(E,F))$ and assume that 
there exists a $B\in\scrL(L^2(\rz^m,F),L^2(\rz^m,E))$ such that 
 $$BA=\Phi-K_1-K_2,$$
where $\Phi=\op(\varphi)\in L^0_\cl(\rz^m,\scrL(E))$ with $\varphi^{(0)}=1$, 
$K_1$ is a compact operator in $L^2(\rz^m,E)$ and $K_2$ induces a continuous operator 
$L^1(\rz^m,E)\to L^2(\rz^m,E)$. Then $a^{(0)}(x_0,\xi_0)$ is injective. 
\end{prop}
}
\begin{proof}
The proof is very similar to the one of Proposition \ref{prop:semi-fred}. {For simplifying notation we again shall assume 
that $E=F$.} Instead of the operator-family $S_\lambda$, 
defined in \eqref{eq:slambda}, we shall now use $S_\lambda\in\scrL(L^2(\rz^{{m}},E))$, $\lambda\ge1$, defined by 
\begin{equation*}
 (S_\lambda u)(x)=\lambda^{{m}/4}e^{i\lambda x\xi_0}u\big(\lambda^{1/2}(x-x_0)\big).
\end{equation*}
Similarly to Section \ref{sec:02.3.1} we can verify that these $S_\lambda$ are isometric isomorphisms and, 
for every $u\in\scrS(\rz^{m},E)$, 
\begin{itemize}
 \item[i$)$] $S^{-1}_\lambda A S_\lambda u\xrightarrow{\lambda\to+\infty}a^{(0)}(x_0,\xi_0)u$ 
  in $L^2(\rz^{m},E)$,\\
  {$S^{-1}_\lambda\Phi S_\lambda u\xrightarrow{\lambda\to+\infty}\varphi^{(0)}(x_0,\xi_0)u=u$ 
  in $L^2(\rz^{m},E)$, } 
 \item[ii$)$] $S_\lambda u\xrightarrow{\lambda\to+\infty} 0$ weakly in $L^2(\rz^{m},E)$, 
 \item[iii$)$] $S_\lambda u\xrightarrow{\lambda\to+\infty} 0$ in $L^1(\rz^{m},E)$.  
\end{itemize}
Now let us choose $u\in\scrS(\rz^{m})$ such that $\|u\|_{L^2(\rz^n)}=1$ and 
{define $u_e$ by $u_e(x)=u(x)e$ with $e\in E$. We obtain 
\begin{align*}
 \|S_\lambda^{-1}\Phi S_\lambda u_e\|_{L^2(\rz^m,E)}=&\|S_\lambda^{-1}((BA+K_1+K_2)S_\lambda u_e\|_{L^2(\rz^m,E)}\\
 \le& \|B\|_{\scrL(L^2(\rz^m,E))}\|S_\lambda^{-1}AS_\lambda u_e\|_{L^2(\rz^m,E)}+\\
     &+ \|K_1 S_\lambda u_e\|_{L^2(\rz^m,E)}  + \|K_2 S_\lambda u_e\|_{L^2(\rz^m,E)}.
\end{align*}
Passing to the limit $\lambda\to+\infty$, using i)--iii) from above, the left-hand side of the latter inequality 
converges to $\|u_e\|_{L^2(\rz^m,E)}=\|e\|_E$, while the right-hand side tends to 
$\|B\|_{\scrL(L^2(\rz^m,E))}\|a^{(0)}(x_0,\xi_0)e\|_E$. We thus derive the estimate 
 $$\|a^{(0)}(x_0,\xi_0)e\|_E\ge \frac{1}{\|B\|_{\scrL(L^2(\rz^m,E))}}\|e\|_E\qquad\forall\;e\in E,$$
which implies the desired injectivity. 
}
\end{proof}

Also the results of Section \ref{sec:02.4} on the spectral invariance extend to the present setting. 
Let us state this explicitly$:$

\begin{thm}\label{thm:analog}
Theorems $\ref{thm:shubin}$, $\ref{thm:spectral}$ and Corollary $\ref{cor:extension}$ remain valid, with obvious 
adaptations, in the framework of bisingular pseudodifferential operators from $L^{\mu,\nu}_\cl(M\times N;E,F)$.   
\end{thm}

%%%%%%%%%%%%%%%%%%%%%%%%%%%%%%%%%%%%%%%%%%%%%%%%%%
%%%%%%%%%%%%%%%%%%%%%%%%%%%%%%%%%%%%%%%%%%%%%%%%%%
\section{Operators of Toeplitz type}\label{sec:04}

Assume we consider a class of operators that act in an associated scale of Sobolev spaces and that in this class 
we can characterize the Fredholm property of an operator by its ellipticity which, by definition, means the invertibility 
of certain principal symbols associated with the operator. It is natural to pose the following problem: Take an 
operator $\wt{A}$ and two projections $P_0,P_1$ in that class of operators $($where projection means that 
$P_j^2=P_j)$, such that the compostion $A=P_1\wt{A}P_0$ makes sense. 
The range spaces of the projections determine closed subspaces of the Sobolev spaces. 
How can we characterize the Fredholm property of $A$, considered 
as an operator acting between these closed subspaces? 

This question has been answered in \cite{Seil}, in a 
quite general context of  ``abstract" pseudodifferential operators. 
We shall apply these results here to the case of bisingular pseudodifferential operators. 
We focus on the case of operators defined on a product $M\times N$ of compact 
manifolds, as described in the preceeding Section \ref{sec:03}; an analogous result also holds true for the class 
of global bisingular operators described in Section \ref{sec:02}. 

Let $E_0$ and $E_1$ be two vector bundles over $M\times N$ and $P_j\in L^{0,0}(M\times N;E_j,E_j)$, $j=0,1$ 
be two projections. The range spaces 
 $$Q^{s,t}(M\times N, E_j;P_j):=P_j\big(Q^{s,t}(M\times N, E_j)\big),\qquad s\in\rz,$$
are closed subspaces of $Q^{s,t}(M\times N, E_j)$. The principal symbols $\sigma^0_0(P_j)$ and $\sigma^0_1(P_j)$, 
see \eqref{eq:hom1} and \eqref{eq:hom2}, are projections when acting as bundle homomorphisms in 
$\pi_M^*\calQ^s(N,E_j)$ and $\pi_N^*\calQ^s(M,E_j)$, respectively. Thus they determine subbundles which 
we shall denote by 
 $$\calQ^s(N,E_j;P_j)\subset \pi_M^*\calQ^s(N,E_j),\qquad 
     \calQ^s(M,E_j;P_j)\subset \pi_N^*\calQ^s(M,E_j).$$
Note that these are bundles on $S^*M$ and $S^*N$, respectively, that generally do not arise as liftings from 
bundles over $M$ and $N$, respectively. 

\begin{thm}\label{thm:proj}
Let $\wt{A}\in L^{\mu,\nu}(M\times N;E_0,E_1)$ and $P_j$ projections as described above. 
For ${A}:=P_1\wt{A}P_0$ the following assertions are equivalent$:$ 
\begin{itemize}
 \item[a$)$] ${A}:Q^{s,t}(M\times N, E_0;P_0)\to Q^{s-\mu,t-\nu}(M\times N, E_1;P_1)$ is a Fredholm operator 
  for some $s\in\rz$. 
 \item[b$)$] The following bundle homomorphisms are isomorphisms$:$
   \begin{align*}  
    \sigma^\mu_0({A})&:\calQ^s(N,E_0;P_0)\lra \calQ^{s-\nu}(N,E_1;P_1),\\ 
    \sigma^\mu_1({A})&:\calQ^s(M,E_0;P_0)\lra \calQ^{s-\mu}(M,E_1;P_1).
   \end{align*} 
\end{itemize}
Moreover, the following two assertions are equivalent$:$
\begin{itemize}
 \item[i$)$] ${A}:Q^{s,t}(M\times N, E_0;P_0)\to Q^{s-\mu,t-\nu}(M\times N, E_1;P_1)$ is invertible  
  for some $s,t\in\rz$. 
 \item[ii$)$] There exists a $\wt{B}\in L^{\mu,\nu}(M\times N;E_1,E_0)$ such that ${A}{B}=P_1$ and 
  ${B}{A}=P_0$ for ${B}:=P_0\wt{B}P_1$. 
\end{itemize}
\end{thm}
\begin{proof}
{
First of all let us observe that we may assume without loss of generality that both bundles $E_0$ and $E_1$ are 
trivial bundles. In fact, due to Swan's theorem, there exists a bundle $E^\prime_0$ over $M$ and such that 
$\calE_0:=E_0\oplus E^\prime_0=M\times N\times\cz^{L_0}$ for some $L_0\in\nz$. 
Similarly, $\calE_1:=E_1\oplus E^\prime_1=M\times N\times\cz^{L_1}$. Now we define the new projections
$\calP_j=\begin{pmatrix}P_j&0\\0&0\end{pmatrix}\in L^{0,0}_\cl(M\times N;\calE_j,\calE_j)$, acting as $P_j$ 
on sections in $E_j$ and as zero on sections in $E_j^\prime$. Similarly, 
we extend $\wt{A}$ to $\wt{\calA}\in L^{\mu,\nu}_\cl(M\times N;\calE_0,\calE_1)$. Then 
$Q^{s,t}(M\times N, \calE_j;\calP_j)=Q^{s,t}(M\times N, E_j;P_j)$ 
and $A$ can be identified with $\calA=\calP_1\wt{\calA}\calP_0$. 
Also the respective principal symbols can be identified with each other.

Next, assuming that the $E_j$ are trivial of fibre-dimension $L_j$,}  
let us justify that we may assume without 
loss of generality that $\mu=\nu=s=t=0$. In fact, let 
$\Lambda^{\sigma,\rho}_j\in L^{\sigma,\rho}(M\times N;E_j,E_j)$, $\sigma,\rho\in\rz$, be invertible with 
$(\Lambda^{\sigma,\rho}_j)^{-1}=\Lambda^{-\sigma,-\rho}_j$.\footnote{{Let 
$\lambda^{\sigma,\rho}=(1-\Delta_M)^{\mu/2}\otimes(1-\Delta_N)^{\nu/2}$ with the Laplacians on $M$ and 
$N$, respectively. Then let $\Lambda^{\sigma,\rho}_j$ be the $(L_j\times L_j)$-diagonal matrix with entries 
$\lambda^{\sigma,\rho}$. }\label{footnote}}
Then the Fredholm property (respectively invertibility$)$ of ${A}$ is equivalent to that of 
 $${A}^\prime:=P_1^\prime \wt{A}^\prime P_0^\prime:
     Q^{0,0}(M\times N, E_0;P_0^\prime)\lra Q^{0,0}(M\times N, E_1;P_1^\prime),$$
where $\wt{A}^\prime:=\Lambda^{s-\mu,t-\nu}_1\wt{A}\Lambda^{-s,-t}_0$ is of zero oder and both  
$P_0^\prime=\Lambda^{s,t}_0P_0\Lambda^{-s,-t}_0$ and 
$P_1^\prime=\Lambda^{s-\mu,t-\nu}_1P_1\Lambda^{\mu-s,\nu-t}_1$ are projections. 

Following \cite{Seil}, let $G:=\{(M\times N;E)\mid E\text{ trivial vector bundle over $M\times N$}\}$, called the set of 
admissible weights, and  
\begin{align*}
 L^\mu(\mathbf{g}):=&L^{\mu,\mu}(M\times N;E_0,E_1),\\ 
     \mathbf{g}=&\big((M\times N;E_0),(M\times N;E_1)\big)\in G\times G
\end{align*}
as well as 
\begin{align*}
 H^s(g):=&Q^{s,s}(M\times N,E),\qquad g=(M\times N;E)\in G.
\end{align*}
Then the equivalence of a$)$ and b$)$ is just Theorem 3.12 of \cite{Seil} $($the assumptions are satisfied due 
to the equivalence of ellipticity and Fredholm property, cf. Theorem \ref{thm:analog} and Section \ref{sec:03.2}$)$, 
while the equivalence of i$)$ and ii$)$ is Theorem 3.9 of \cite{Seil}. 
\end{proof}

{
\subsection{Order reductions}\label{sec:04.1}

In this section we shall show the existence of bisingular order reductions on a product of two closed manifolds. 
We shall need the following lemma: 

\begin{lem}\label{lem:positive} 
Let $\mu>0$  and $A\in L^\mu_\cl(M,\cz^L)$ be elliptic, symmetric and have scalar principal symbol. 
Moreover, assume that $A$ is positive, i.e., 
 $$(Au,u)_{L^2(M,\cz^L)}>0\qquad\forall\;0\not=u\in\scrC^\infty(M,\cz^L).$$
Let $P\in L^0_\cl(M,\cz^L)$ be an orthogonal projection. Then 
 $$A_P:=PAP+(1-P)A(1-P)\in L^\mu_\cl(M,\cz^L)$$
is invertible with inverse belonging to $L^{-\mu}_\cl(M,\cz^L)$. 
\end{lem}
\begin{proof}
Since $A$ has scalar principal symbol, $A_P$ has the same principal symbol as $A$, hence is elliptic. 
Since $P$ is orthogonal, $A_P$ is also positive. It remains to observe that the spectrum of elliptic operators 
of positive order consists of isolated eigenvalues only. Due to the positivity, $0$ is not an eigenvalue of $A_P$. 
\end{proof}

\begin{thm}\label{thm:ord-red}
Let $\mu,\nu\in\rz$ and $E$ be a Hermitian vector bundle over $M\times N$. Then there exist operators 
$A\in L^{\mu,\nu}_\cl(M\times N;E,E)$ and $B\in L^{-\mu,-\nu}_\cl(M\times N;E,E)$ such that 
$AB=1$ and $BA=1$. 
\end{thm}
}

{
Observe that it is sufficient to show this theorem in case $\mu,\nu>0$. 
In fact, given arbitrary $\mu,\nu$ choose $\mu_0,\nu_0>0$ such that $\mu_1:=\mu+\mu_0>0$ and 
$\nu_1:=\nu+\nu_0>0$. Then choose $A_0\in L^{\mu_0,\nu_0}_\cl(M\times N;E,E)$ and 
$A_1\in L^{\mu_1,\nu_1}_\cl(M\times N;E,E)$ with corresponding inverses $B_0$ and $B_1$. Then 
$A:=B_0A_1\in L^{\mu,\nu}_\cl(M\times N;E,E)$ and $B:=B_1A_0\in L^{-\mu,-\nu}_\cl(M\times N;E,E)$ are as 
desired. 

\begin{proof}[Proof of Theorem \ref{thm:ord-red}]
Let $\mu,\nu>0$. As described in the beginning of the proof of Theorem \ref{thm:proj}, we find a bundle $E^\prime$ 
over $M\times N$ such that $E\oplus E^\prime=M\times N\times \cz^{L}$ with an orthogonal direct sum. 
Let $P$ denote the orthogonal projection onto $E$ along $E^\prime$; we consider $P$ as an element of 
$L^{0,0}_\cl(M\times N;\cz^{L},\cz^{L})$. Then we have the identification 
 $$L^{\mu,\nu}_\cl(M\times N;E,E)
     =\Big\{P\wt{A}P \mid \wt{A} \in L^{\mu,\nu}_\cl(M\times N;\cz^{L},\cz^{L})\Big\},$$
where $P\wt{A}P$ is considered as a map in $($the identified spaces$)$ 
 $$Q^{s,t}(M\times N, E)=Q^{s,t}(M\times N, \cz^{L};P).$$  
Now let $\Lambda=\Lambda^{\mu,\nu}\in L^{\mu,\nu}_\cl(M\times N;\cz^L,\cz^L)$ be as described in 
footnote \ref{footnote}. In particular, $\Lambda$ is elliptic, symmetric and is positive, i.e.,   
 $$(\Lambda u,u)_{L^2(M\times N,\cz^L)}>0\qquad \forall\;0\not=u\in\scrC^\infty(M\times N,\cz^L).$$
Let $\Lambda_P:=P\Lambda P+(1-P)\Lambda(1-P)$. By Lemma \ref{lem:positive} $($applied pointwise/fibrewise to 
the principal symbols $\sigma_1^\mu(\Lambda_P)$ and $\sigma_2^\nu(\Lambda_P)$ of $\Lambda_P)$, 
one sees that $\Lambda_P\in L^{\mu,\nu}_\cl(M\times N;\cz^L,\cz^L)$ is elliptic. 
Moreover, $\Lambda_P$ is symmetric and positive.  
Since the spectrum of elliptic bisingular pseudodifferential operators of positive order$($s$)$ consists of isolated, 
positive eigenvalues $($due to the compact embedding of Sobolev spaces of positive order$($s$)$ into $L^2)$, 
and due to the spectral invariance of bisingular operators, we conclude that $\Lambda_P$ 
is invertible with inverse in $L^{-\mu,-\nu}_\cl(M\times N;\cz^L,\cz^L)$. Then 
$A:=P\Lambda_PP=P\Lambda P$ induces isomorphisms $Q^{s,t}(M\times N,\cz^L;P)\to 
Q^{s-\mu,t-\nu}(M\times N,\cz^L;P)$, i.e., $Q^{s,t}(M\times N,E)\to Q^{s-\mu,t-\nu}(M\times N,E)$. 
Now, due to Theorem \ref{thm:proj}, there exists a $B=P\wt{B}P$ with 
$\wt{B}\in L^{-\mu,-\nu}_\cl(M\times N;\cz^{L},\cz^{L})$ 
such that $AB=BA=P$, hence $AB=BA=1$ on any $Q^{s,t}(M\times N, \cz^{L};P)=Q^{s,t}(M\times N, E)$. 
\end{proof}
}
%%%%%%%%%%%%%%%%%%%%%%%%%%%%%%%%%%%%%%%%%%%%%%%%%%
%%%%%%%%%%%%%%%%%%%%%%%%%%%%%%%%%%%%%%%%%%%%%%%%%%
\section{Appendix: A remark on vector bundles over product spaces}\label{sec:05}

Let $E$ be a vector bundle over $M\times N$, the product of two smooth closed manifolds. 
For every $m\in M$ we define an embedding of $N$ into $M\times N$ by  
 $$\iota_m:N\to M\times N,\quad n\mapsto (m,n)$$
and we denote by $E(m):=\iota_m^* E$ be the corresponding pull-back of $E$ to $N$. 

\begin{prop}\label{prop:triv}
For every $m\in M$ exists an open neighborhood $U\subset M$ such that 
$E|_{U\times N}\cong U\times E({m})$ $($diffeomorphism between smooth manifolds$)$. 
\end{prop}
\begin{proof}
By Swan's theorem we may assume that $E$ is a subbundle of $M\times N\times\cz^N$ for some $N\in\nz$. 
Hence there exists a function $p\in\scrC^\infty\big(M\times N,\scrL(\cz^N)\big)$ taking values in the 
projections of $\cz^N$ and such that 
 $$E_{(m,n)}=\big\{\big(m,n,p(m,n)v\big)\mid v\in\cz^N\big\},\quad 
     E(m)_n=\big\{\big(n,p(m,n)v\big)\mid v\in\cz^N\big\}$$
are the fibres of $E$ over $(m,n)$ and of $E(m)$ over $n$, respectively. 
Now let $m_0\in M$ be fixed. Define $\varphi\in\scrC^\infty\big(M\times N,\scrL(\cz^N)\big)$ by 
 $$\varphi(m,n)=p(m_0,n)+(1-p)(m,n).$$
Since $\varphi(m_0,n)=1$ for every $n$ and since $N$ is compact, we find an open neighborhood 
$U_0$ of $m_0$ such that $\varphi(m,n)\in\scrL(\cz^N)$ is an isomorphism for every $(m,n)\in U_0\times N$. 
In particular, $\varphi$ induces a bundle isomorphism $\Phi$ in $U_0\times N\times\cz^N$. Moreover, 
 $$\Phi(E_{(m,n)})=\{m\}\times E(m_0)_n,\qquad (m,n)\in U_0\times N.$$
In fact, since both sides have the same dimension, this follows if the left-hand side is a subset of the 
right-hand side. However, this is true, since $\varphi(m,n)p(m.n)v=p(m_0,n)p(m,n)v\in\mathrm{im}\,p(m_0,n)$ 
for every $v\in\cz^N$. In other terms, we have verified that $\Phi:E|_{U_0\times N}\to U_0\times E({m_0})$ 
diffeomorphically. 
\end{proof}

\begin{cor}\label{cor:fibre}
Let $M$ be connected and $m_0\in M$ be fixed. Then$:$
\begin{itemize}
 \item[a$)$] $E(m)$ is isomorphic to $E(m_0)$ for every $m\in M$. 
 \item[b$)$] $E$ is a fibre bundle over $M$ with typical fibre $E(m_0)$. 
\end{itemize} 
\end{cor}
\begin{proof}
For a$)$ denote by $V$ the set of all $m\in M$ such that $E(m)\cong E(m_0)$. By Proposition \ref{prop:triv} both 
$V$ and $M\setminus V$ are open subsets of $M$. Since $m_0\in M$ and $M$ is connected, $M\setminus V$ must 
be empty, hence $V=M$. Clearly, b$)$ follows from a$)$ and Proposition \ref{prop:triv}. 
\end{proof}

In the following let $Q^s(N,F)$ denote the standard $L^2$-Sobolev space of order $s$ of sections in the vector 
bundle $F$ over $N$. This is a separable, infinite dimensional Hilbert space.     

\begin{cor}
Let $m_0\in M$ be fixed $($and $M$ not necessarily connected$)$. Then 
 $$\calQ^s(N,E):=\mathop{\mbox{\Large$\cup$}}_{m\in M}\{m\}\times Q^s(N,E(m))$$
is a Hilbert space bundle over $M$ with typical fibre $Q^s(N,E(m_0))$. 
\end{cor}
\begin{proof}
Let $M_0,\ldots,M_k$ be the connected components of $M$ and fix points $m_i\in M_i$. 
Corollary \ref{cor:fibre} implies that $\calQ^s(N,E)|_{M_i}$ is a bundle over $M_i$ with typical fibre  
$Q^s(N,E(m_i))$. It remains to observe that any $Q^s(N,E(m_i))$ is isomorphic to $Q^s(N,E(m_0))$, 
since all these spaces are isomorphic to $\ell^2(\nz)$, for example.  
\end{proof}

% ------------------------------------------------------------------------
\end{document}